\documentclass[11pt]{article}
\usepackage{amssymb} %Sonderzeichen
\usepackage{amsthm} %theoremstyle
\usepackage{amsmath} %numberwithin
\usepackage{float} %Table position
\usepackage{fancyhdr} %Eigener Header
\usepackage[english]{babel}
\usepackage[utf8]{inputenc}
\usepackage[T1]{fontenc}
\usepackage[margin=1.5in]{geometry} %Rand
\usepackage{wasysym}

\usepackage[osf,sc]{mathpazo}
\usepackage{textcomp}
\usepackage{microtype}
\usepackage{graphicx, color}
\usepackage{amsfonts}
\usepackage{booktabs}
\usepackage{pifont}
\usepackage{tikz}
\usepackage{graphicx}
\usepackage{float}
\usepackage{hyperref}
\usepackage[all]{xy}

\begin{document}
\title{Euler and the Duplication Formula for the Gamma-Function}
\author{Alexander Aycock, Johannes-Gutenberg University Mainz \\ Staudinger Weg 9, 55128 Mainz \\ \url{aaycock@students.uni-mainz.de} }
\date{ }
\maketitle

\begin{abstract}
We show how the formulas in Euler's paper {}"Variae considerationes circa series hypergeometricas"{} \cite{E661} imply Legendre's duplication formula for the $\Gamma$-function. This paper can be seen as an Addendum to \cite{Ay21}.
\end{abstract}

\section{Introduction}
\label{sec: Introduction}

In \cite{Ay21}, we focused on a function defined by Euler in \cite{E661} as:

\begin{equation}
\label{eq: Gamma General}
    \Gamma_E(x):= a\cdot (a+b) \cdot (a+2b) \cdot (a+3b) \cdot \cdots \cdot (a+(x-1)b) \quad \text{for} \quad a,b >0,
\end{equation}
which we showed to be continueable to non-integer values of $x$ via the expression:

\begin{equation}
\label{eq: General Gamma via Gamma}
    \Gamma_E(x)= \dfrac{b^x}{\Gamma \left(\frac{a}{b}\right)}\cdot \Gamma \left(x+\dfrac{a}{b}\right).
\end{equation}
Here, $\Gamma(x)$ means the ordinary $\Gamma$-function defined as:

\begin{equation}
    \label{eq: Gamma Function}
    \Gamma(x) := \int\limits_{0}^{\infty}e^{-t}t^{x-1}dt \quad \text{for} \quad \operatorname{Re}(x)>0.
\end{equation}
Equation (\ref{eq: Gamma General}) enabled us to determine the constant $A$ in the asymptotic expansion for the function $\Gamma_E$ found by Euler via the Euler-Maclaurin summation formula. The asymptotic expansions reads:

\begin{equation}
    \Gamma_E(x) \sim A \cdot e^{-x}\cdot (a-b+bx)^{\frac{a}{b}+x-\frac{1}{2}} \quad \text{for} \quad x \rightarrow \infty.
\end{equation}
We found the constant $A$ to be

\begin{equation}
    \label{eq: Constant A}
    A= \dfrac{\sqrt{2\pi}}{\Gamma\left(\frac{a}{b}\right)}\cdot e^{1-\frac{a}{b}}\cdot b^{\frac{1}{2}-\frac{a}{b}}.
\end{equation}
In this paper, we intend to use this result and more of Euler's formulas from the same paper to show that they imply the Legendre duplication formula for the $\Gamma$-function, i.e., the relation

\begin{equation}
    \label{eq: Legendre Duplication Formula}
    \Gamma(x) = \dfrac{2^{x-1}}{\sqrt{\pi}}\cdot \Gamma\left(\dfrac{x}{2}\right)\cdot \Gamma \left(\dfrac{x}{2}+\dfrac{1}{2}\right).
\end{equation}

\section{Euler's other Functions}
\label{sec: Euler's other Functions}

\subsection{Euler's Definition}
\label{subsec: Euler's Definition}

Aside from the function $\Gamma_E$, in his paper \cite{E661}, Euler introduced two other related functions:

\begin{equation}
\label{eq: Euler's other Functions}
        \renewcommand{\arraystretch}{1,5}
\setlength{\arraycolsep}{0.0mm}
\begin{array}{llllll}
    \Delta(x) &~=~& a \cdot (a+2b)\cdot (a+4b)\cdot (a+6b)\cdot \cdots \cdot (a+(2x-2)b), \\
    \Theta(x) &~=~& (a+b)\cdot (a+3b) \cdot (a+5b) \cdot \cdots \cdot (a+(2x-1)b).
\end{array}
\end{equation}
As it was the case for $\Gamma_{E}$ (equation (\ref{eq: Gamma General})), Euler's definition is only valid for integer values of $x$, but by using the ideas from \cite{Ay21}, we could extend the definition to real numbers. 

\subsection{Asymptotic Expansions of these Functions}
\label{subsec: Asymptotic Expansions of these Functions}

Furthermore, Euler also found asymptotic expansions for his functions $\Delta$ and $\Theta$. They are:

\begin{equation}
\label{eq: Delta}
        \renewcommand{\arraystretch}{1,5}
\setlength{\arraycolsep}{0.0mm}
\begin{array}{llllll}
     \Delta(x) &~\sim ~& B\cdot e^{-x} \cdot (a-2b+2bx)^{\frac{a}{2b}+x-\frac{1}{2}} \\
     \Theta(x) &~\sim ~& C \cdot e^{-x} \cdot (a-b+2bx)^{\frac{a}{2b}+x},
\end{array}
\end{equation}
where $B$ and $C$ are constants resulting from the application of the Euler-Maclaurin summation formula and the asymptotic expansions are valid for $x \rightarrow \infty$.

\subsection{Relation among the Constants}
\label{subsubsec: Relation among the Constants}

Euler was not able to find any of the constants $A$, $B$ and $C$. But, using the general relations among his functions $\Gamma_E$, $\Delta$ and $\Theta$ and the respective corresponding asymptotic expansions, he found the following relations:

\begin{equation}
    \label{eq: Relation among the Constants A}
   A= \dfrac{B \cdot C}{\sqrt{e}} 
\end{equation}
and

\begin{equation}
     \label{eq: Relation among the Constants B}
     B = C\cdot k \cdot \sqrt{e}
\end{equation}
with $ k=\Delta\left(\frac{1}{2}\right)$.
As we will show in the next section, these relations imply the Legendre duplication formula (equation (\ref{eq: Legendre Duplication Formula})).

\section{Derivation of the Legendre Duplication Formula from Euler's Formulas}
\label{sec: Derivation of the Legendre Duplication Formula from Euler's Formulas}

As Euler remarked himself in \cite{E661}, equations (\ref{eq: Relation among the Constants A}) and (\ref{eq: Relation among the Constants B}) tell us that we only need to find one of the constants $A$, $B$ and $C$ such that we can calculate the remaining two from the first. Since we discovered the value $A$ (equation (\ref{eq: Constant A})), we could do precisely that. But for our task at hand, we need to find the value of $k$ first.

\subsection{Evaluation of the Constant $k$}
\label{subsec: Evaluation of the Constant $k$}

To evaluate $k=\Delta\left(\frac{1}{2}\right)$, we note that we just have to make the substitution $b \mapsto 2b$ in equation (\ref{eq: Gamma General}) such that the expression for $\Gamma_E$ goes over into the expression for $\Delta$ (equation (\ref{eq: Delta})) in equation (\ref{eq: Euler's other Functions}). Making the same substitution in equation (\ref{eq: General Gamma via Gamma}), we arrive the the following expression for $\Delta(x)$:

\begin{equation*}
    \Delta(x) = \dfrac{(2b)^x}{\Gamma\left(\frac{a}{2b}\right)}\cdot \Gamma \left(x+\dfrac{a}{2b}\right).
\end{equation*}
Therefore, for $x=\frac{1}{2}$

\begin{equation}
\label{eq: Value k}
    k = \Delta \left(\dfrac{1}{2}\right)=  \dfrac{(2b)^{\frac{1}{2}}}{\Gamma\left(\frac{a}{2b}\right)}\cdot \Gamma \left(\dfrac{1}{2}+\dfrac{a}{2b}\right).
\end{equation}

\subsection{The Legendre Duplication Formula}
\label{subsec: The Legendre Duplication Formula}

Having found $k$, let us use equations (\ref{eq: Relation among the Constants A}) and (\ref{eq: Relation among the Constants B}) to find the Legendre duplication formula (equation (\ref{eq:  Legendre Duplication Formula})). Substituting the value for $C$ in (\ref{eq: Relation among the Constants B}) in for the value of $C$ in (\ref{eq: Relation among the Constants A}), we arrive at this equation:

\begin{equation}
\label{eq: Relation AB}
    A= \dfrac{B^2}{\Delta \left(\frac{1}{2}\right)}e^{-1}.
\end{equation}
Next, we note that since $\Delta(x)$ is obtained from $\Gamma_{E}(x)$ by the substitution $b \mapsto 2b$, the value of the constant $B$ is obtained in the same way from $A$ and reads:

\begin{equation}
    \label{eq: Constant B}
    B = \dfrac{\sqrt{2\pi}}{\Gamma \left(\frac{a}{2b}\right)}\cdot (2b)^{\frac{1}{2}-\frac{a}{2b}}\cdot e^{1-\frac{a}{2b}}.
\end{equation}
Thus, substituting the respective values for $A$ (equation (\ref{eq: Constant A})), $B$ (equation (\ref{eq: Constant B})) and $k$ (equation (\ref{eq: Value k})), equation (\ref{eq: Relation AB}) becomes:

\begin{equation*}
    \dfrac{\sqrt{2\pi}}{\Gamma\left(\frac{a}{b}\right)}\cdot e^{1-\frac{a}{b}}\cdot b^{\frac{1}{2}-\frac{a}{b}} = \dfrac{\left(\frac{\sqrt{2\pi}}{\Gamma \left(\frac{a}{2b}\right)}\cdot (2b)^{\frac{1}{2}-\frac{a}{2b}}\cdot e^{1-\frac{a}{2b}}\right)^2}{\frac{(2b)^{\frac{1}{2}}}{\Gamma\left(\frac{a}{2b}\right)}\cdot \Gamma \left(\frac{1}{2}+\frac{a}{2b}\right)} \cdot e^{-1}.
\end{equation*}
Most terms cancel each other and after this equation simplifies to:

\begin{equation*}
    \dfrac{1}{\Gamma \left(\frac{a}{b}\right)}= \dfrac{\sqrt{2\pi}\cdot 2^{\frac{1}{2}-\frac{a}{b}}}{\Gamma \left(\frac{a}{2b}\right)\cdot \Gamma \left(\frac{1}{2}+\frac{a}{2b}\right)}.
\end{equation*}
Finally, writing $x$ instead of $\frac{a}{b}$ and solving this equation for $\Gamma(x)$, after a little simplification, we arrive at the relation:

\begin{equation*}
    \Gamma(x) = \dfrac{2^{x-1}}{\sqrt{\pi}}\cdot \Gamma \left(\dfrac{x}{2}\right)\cdot \Gamma \left(\dfrac{x+1}{2}\right),
\end{equation*}
which is the Legendre duplication formula for the $\Gamma$-function (equation (\ref{eq: Legendre Duplication Formula})), as we wanted to show.

\section{Conclusion}
\label{sec: Conclusion}

In this note we showed that Legendre's duplication formula, i.e., equation (\ref{eq: Legendre Duplication Formula}) follows from Euler's formulas found in his paper \cite{E661}. Indeed, the Legendre duplication formula could also have been shown by Euler himself, if he had set this task for himself, as we argued in more detail in \cite{Ay21}.
Furthermore, Euler's ideas that we explained in this  and the before-mentioned paper, can be generalized to show the multiplication formula for the $\Gamma$-function, i.e, the formula

\begin{equation*}
    \label{Multiplication Formula}
    \Gamma(x) = \sqrt{\frac{n}{(2\pi)^{n-1}}}\cdot n^{x-1} \cdot \Gamma\left(\dfrac{x}{n}\right)\Gamma\left(\dfrac{x+1}{n}\right)\Gamma\left(\dfrac{x+2}{n}\right) \cdot \cdots \cdot \Gamma\left(\dfrac{x+n-1}{n}\right).
\end{equation*}
This formula is attributed to Gauss who stated and proved it in \cite{Ga28}. But it was given  by Euler (in different form, expressed via Beta functions) in \cite{E421}, as we demonstrated in \cite{Ay19}.

\end{document}